\documentclass[letterpaper]{article}
 
\title{The topology of competitively constructed graphs}
\date{January 15, 2013}
\author{Alan Frieze\thanks{Research supported in part by NSF grant ccf1013110},  Wesley Pegden}

\usepackage{amsmath,amsthm,amssymb,amscd,latexsym,pstricks,pst-node,pst-tree}
\usepackage{subfigure}
\usepackage{url}
\usepackage{cite}

\renewcommand{\ldots}{\dots}

\def\cT{{\cal T}}

\newcommand{\tref}[1]{Theorem~\ref{t.#1}}
\newcommand{\oref}[1]{Observation~\ref{o.#1}}

\newcommand{\lref}[1]{Lemma~\ref{l.#1}}
\newcommand{\fref}[1]{Figure~\ref{f.#1}}

\newcommand{\eref}[1]{(\ref{e.#1})}

\newcommand{\sbs}{\subset}

\newtheorem{theorem}{Theorem}[section]

\newtheorem{lemma}[theorem]{Lemma}

\newtheorem{observ}[theorem]{Observation}

\newcommand{\comments}[1]{}

{
\theoremstyle{definition}

\newtheorem{q}{Question
}

{
\theoremstyle{remark}

}

\newcommand{\stm}{\setminus}

\newcommand{\df}{\mathrm{def}}

\begin{document}
\maketitle

\begin{abstract}
We consider a simple game, the \emph{$k$}-regular graph game, in which players take turns adding edges to an initially empty graph subject to the constraint that the degrees of vertices cannot exceed $k$.  We show a sharp topological threshold for this game: for the case $k=3$ a player can ensure the resulting graph is planar, while for the case $k=4$, a player can force the appearance of arbitrarily large clique minors.
\end{abstract}

\section{Introduction}
In some sense, restricting one's attention to 3-regular graphs is not a topological constraint at all, in the sense that connected 3-regular graphs can require arbitrarily complex surfaces to be embedded in, or, say, contain arbitrarily large clique minors.  In particular, from a topological point of view, vertices of degree 3 are essentially different from vertices of degree 2.  One might then expect the presence of degrees greater than 3 to lead to a similar increase in topological trouble.  For example, out of the list of 103 forbidden topological minors for embeddability in the projective plane, only 6 are required to ensure embeddability of cubic graphs in the projective plane\cite{glover1975cubic,archdeacon1981kuratowski}.

We show another kind of topological threshold between degree-3 and degree-4 graph vertices.  Consider a game (the \emph{$k$-regular graph game}) in which two players take turns adding edges to an initially empty graph.  Players are allowed to add edges only between pairs of vertices which were previously nonadjacent and of degree $\leq k-1$.  The game ends when this is no longer possible.  In particular, the degree of every vertex in the resulting graph will be exactly $k$, with at most $k$ exceptions.  

\begin{theorem}
\label{t.3game}
Regardless of who has the first move, a player in the 3-regular graph game has a strategy to ensure that the resulting graph is planar.
\end{theorem}
On the other hand, for the analogous $4$-regular graph game, we have:
\begin{theorem}
\label{t.4game}
For any $\ell$ and sufficiently large $n$, and regardless of who has the first move, a player in the 4-regular graph game on $n$ vertices has a strategy to ensure that the resulting graph has a $K_\ell$ minor.
\end{theorem}
\noindent Thus there is no surface $S$ for which a player of the 4-regular graph game can ensure that the connected components of the result has a drawing on $S$.

\smallskip 
Note that the moves of the two players in this game are equivalent, unlike Maker-Breaker games (see \cite{Beck}).  In particular, with a symmetric ``normal'' win condition---say, the first player in the 6-regular game to break planarity loses---this would be an impartial game, subject to the Sprague-Grundy theorem \cite{sprague,grundy}.
\section{Proofs}
We begin by proving \tref{3game}.  Call the player with the goal of making the graph planar the \emph{planar player}; his opponent is the \emph{nonplanar player}.  At any stage of the game, the \emph{deficit} $\df(v)$ of a vertex $v$ refers to the difference between the current degree and the maximum allowable degree.  Thus, in the 3-regular graph game, every vertex begins with deficit 3.  The deficit of a set of vertices is the sum of their deficits.

We inductively claim that the planar player can maintain that at any stage, at the end of his move, each connected component $C$ of $G$ can be drawn in the plane such that its positive deficit vertices all lie on its unbounded face, and also that each $C$ is one of the following \emph{Types}:

\psset{unit=.5,yunit=-1,linewidth=.5pt,dotsize=5pt}

\begin{figure}[t]
\centering
\begin{pspicture}(-7,-4)(10,1)
\SpecialCoor

\pscircle(-6,-3){1}
\rput(-6,-3){
\psdot(1;90)
\psdot(1;-30)
\psdot(1;210)
}

\pscircle(-3,-3){1}
\psdot(-3,-2)
\psdot(-3,-4)

\pscurve(-2,-3)(-1.8,-2.8)(-1.3,-3.3)(-1,-3)

\pscircle(0,-3){1}
\psdot(0,-2)
\psdot(0,-4)

\pscircle(3,-3){1}
\psdot(3,-2)
\psdot(3,-4)

\pscurve(4,-3)(4.2,-2.8)(4.7,-3.3)(5,-3)

\pscircle(6,-3){1}
\psdot(6,-4)

\pscurve(7,-3)(7.2,-2.8)(7.7,-3.3)(8,-3)

\pscircle(9,-3){1}
\psdot(9,-4)
\psdot(9,-2)

\pscircle(-3,0){1}
\psdot(-3,1)
\psdot(-3,-1)

\pscurve(-2,0)(-1.8,.2)(-1.3,-.3)(-1,0)

\pscircle(0,0){1}
\psdot(0,-1)

\pscurve(1,0)(1.2,.2)(1.7,-.3)(2,0)

\pscircle(3,0){1}
\psdot(3,-1)

\pscircle(6,0){1}
\psdot(6,1)
\psdot(6,-1)

\pscurve(4,0)(4.2,.2)(4.7,-.3)(5,0)

\end{pspicture}
\caption{\label{f.types} Components of Types 1, 2, \ref{small3} and \ref{big3}, respectively}
\end{figure}
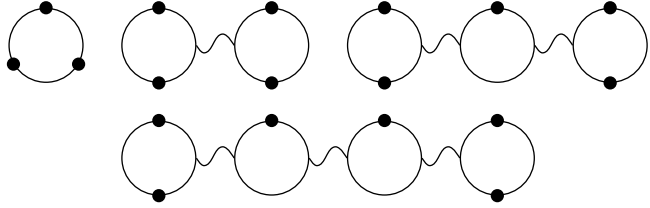

\begin{enumerate}
\item \label{t1} $C$ has deficit $\leq 3$;
\item \label{t2} $C$ has a bridge $e$, such that the vertex-sets of the connected components $C_1,C_2$ of $C\stm e$ (the \emph{sides} of $C$) each have deficit exactly 2; or
\item \label{t3} 
\begin{enumerate}
\item \label{small3} $C$ has bridges $e_1$, $e_2$,  such that the component $C_0$ of $C\stm e_1$ which does not contain $e_2$ has deficit 2, the component $C_2$ of $C\stm e_2$ which does not contain $e_1$ has deficit 2, and the component $C_1$ of $C\stm \{e_1,e_2\}$ which is not $C_0$ or $C_2$ has deficit 1.
\item \label{big3} $C$ has bridges $e_1$, $e_2$, $e_3$, where $e_1$ and $e_3$ are in distinct components of $C\stm e_2$,  such that the component $C_0$ of $C\stm e_1$ which does not contain $e_2,e_3$ has deficit 2, the component $C_3$ of $C\stm e_3$ which does not contain $e_1,e_2$ has deficit 2, and the components $C_1,C_2$ of $C\stm \{e_1,e_2,e_3\}$ which are not $C_0$ or $C_3$ have deficit 1.
\end{enumerate}
\end{enumerate}
Moreover, the planar player will ensure that there is at most one Type
3 component after each of his turns (which is why we consider Types \ref{small3} and \ref{big3} to be variants of a single type).
At any stage, if all components are of one of these types and there is at most one
component of Type 3, we say the graph satisfies condition $\cT$. 

\fref{types} shows schematics of these component types.  Note that if a component $C_i$ has deficit 2, it has either two vertices of deficit 1, or a single vertex of deficit 2; similarly, when a component has deficit 3, it may be distributed among vertices in 3 different ways.  Our argument will not be sensitive to these distinct cases.  In particular, the drawings in \fref{types} and \fref{moves} use circles to denote components after the deletion of the relevant bridges, and dots to denote units of deficit within these components, but we do not assume that distinct dots drawn on individual circles correspond to distinct vertices, and it is thus not valid to assume, for example, that two components with deficit $\geq 2$ can be joined by two distinct edges.  By the same token, we may not assume that it is legal to draw an edge between the two sides of a Type 2 component (it may be a single edge, for example).

\smallskip

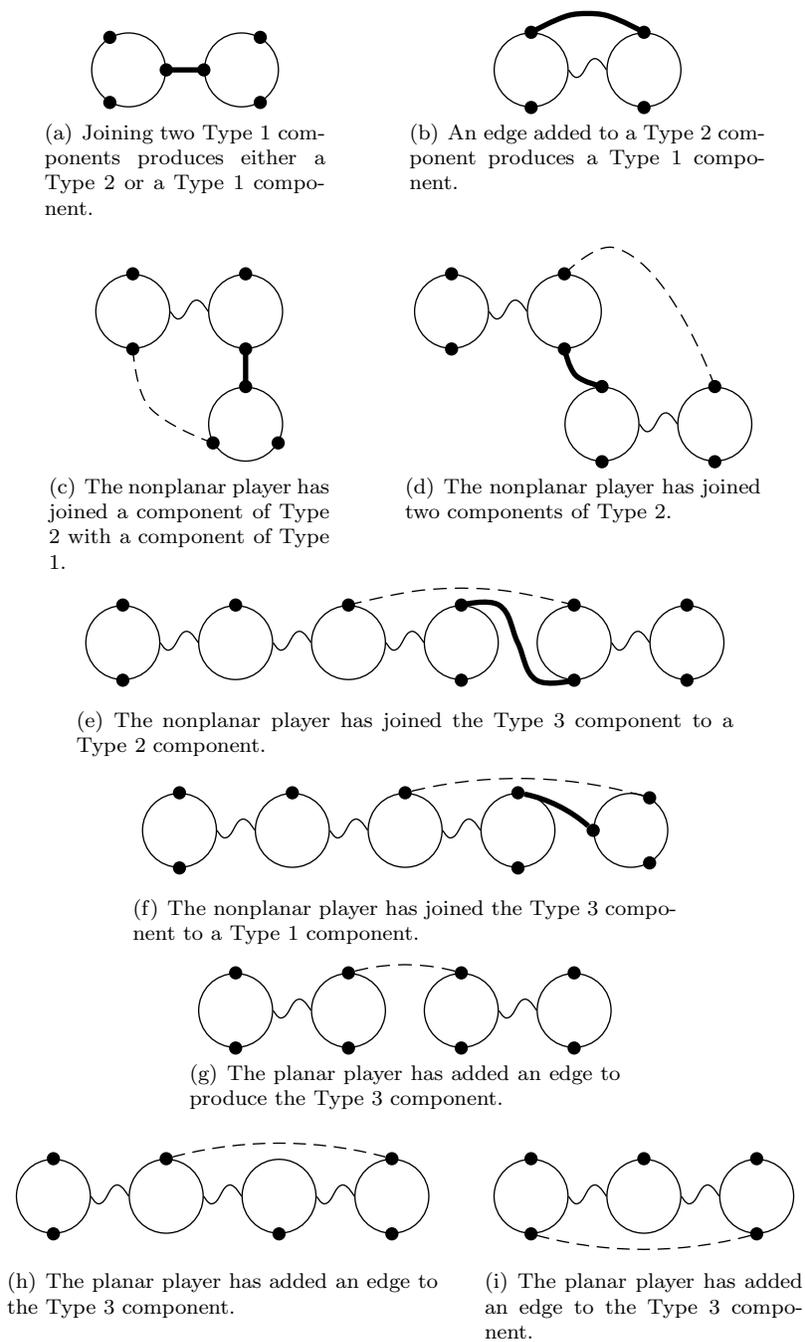
\begin{figure}[p]

\newcommand{\vx}[2]{\cnode*(#1){2.5pt}{#2}}
\newcommand{\Dedg}[2]{\ncline[linestyle=dashed]{#1}{#2}}
\newcommand{\Dcedg}[2]{\ncarc[linestyle=dashed,arcangleA=15,arcangleB=15]{#1}{#2}}
\newcommand{\Bedg}[2]{\ncline[linewidth=2pt]{#1}{#2}}
\newcommand{\Bcedg}[2]{\ncarc[linewidth=2pt,arcangleA=15,arcangleB=15]{#1}{#2}}

\begin{center}

\subfigure[\label{f.type1type1}Joining two Type 1 components produces either a Type 2 or a Type 1 component.]{
\begin{pspicture}(-2,-1)(5,1)
\SpecialCoor

\pscircle(0,0){1}
\rput(0,0){
\psdot(1;0)
\psdot(1;120)
\psdot(1;240)
}

\psline[linewidth=2pt](1,0)(2,0)

\pscircle(3,0){1}
\rput(3,0){
\psdot(1;180)
\psdot(1;60)
\psdot(1;300)
}
\end{pspicture}
}\hspace{1cm} %
\subfigure[\label{f.type2self}An edge added to a Type 2 component produces a Type 1 component.]{
\begin{pspicture}(-3,-1)(6,1)
\SpecialCoor

\pscircle(0,0){1}
\psdot(0,1)
\psdot(0,-1)

\pscurve(1,0)(1.2,.2)(1.7,-.3)(2,0)

\pscircle(3,0){1}
\psdot(3,1)
\psdot(3,-1)

\pscurve[linewidth=2pt](3,-1)(2,-1.45)(1.5,-1.5)(1,-1.45)(0,-1)
\end{pspicture}
}

\subfigure[\label{f.type2type1}The nonplanar player has joined a component of Type 2 with a component of Type 1.]{
\begin{pspicture}(-2,-1.8)(5,4)
\SpecialCoor

\pscircle(0,0){1}
\psdot(0,1)
\psdot(0,-1)

\pscurve(1,0)(1.2,.2)(1.7,-.3)(2,0)

\pscircle(3,0){1}
\psdot(3,1)
\psdot(3,-1)

\pscircle(3,3){1}
\rput(3,3){
\psdot(1;90)
\psdot(1;210)
\psdot(1;330)
}
\psline[linewidth=2pt](3,1)(3,2)

\pscurve[linestyle=dashed](0,1)(.4,2.5)(2.133,3.5)

\end{pspicture}
}\hspace{1cm}%
\subfigure[\label{f.type2type2}The nonplanar player has joined two components of Type 2.]{
\begin{pspicture}(-1,-1.8)(8,4)
\SpecialCoor

\pscircle(0,0){1}
\psdot(0,1)
\psdot(0,-1)

\pscurve(1,0)(1.2,.2)(1.7,-.3)(2,0)

\pscircle(3,0){1}
\psdot(3,1)
\psdot(3,-1)

\pscircle(4,3){1}
\psdot(4,4)
\psdot(4,2)

\pscurve(5,3)(5.2,3.2)(5.7,2.7)(6,3)

\pscircle(7,3){1}
\psdot(7,4)
\psdot(7,2)

\pscurve[linewidth=2pt](3,1)(3.3,1.7)(4,2)

\pscurve[linestyle=dashed](3,-1)(4.3,-1.7)(7,2)

\end{pspicture}
}

\subfigure[\label{f.type3type2}The nonplanar player has joined the Type 3 component to a Type 2 component.]{
\begin{pspicture}(-4,-1.2)(13,1.4)
\SpecialCoor

\pscircle(-3,0){1}
\psdot(-3,1)
\psdot(-3,-1)

\pscurve(-2,0)(-1.8,.2)(-1.3,-.3)(-1,0)

\pscircle(0,0){1}
\psdot(0,-1)

\pscurve(1,0)(1.2,.2)(1.7,-.3)(2,0)

\pscircle(3,0){1}
\psdot(3,-1)

\pscircle(6,0){1}
\psdot(6,1)
\psdot(6,-1)

\pscurve(4,0)(4.2,.2)(4.7,-.3)(5,0)

\pscircle(9,0){1}
\psdot(9,1)
\psdot(9,-1)

\pscurve(10,0)(10.2,.2)(10.7,-.3)(11,0)

\pscircle(12,0){1}
\psdot(12,1)
\psdot(12,-1)

\pscurve[linewidth=2pt](6,-1)(7,-1)(7.5,0)(8,1)(9,1)

\pscurve[linestyle=dashed](3,-1)(5,-1.4)(7,-1.4)(9,-1)

\end{pspicture}
}

\subfigure[\label{f.type3type1}The nonplanar player has joined the Type 3 component to a Type 1 component.]{
\begin{pspicture}(-4,-1.2)(10,1.4)
\SpecialCoor

\pscircle(-3,0){1}
\psdot(-3,1)
\psdot(-3,-1)

\pscurve(-2,0)(-1.8,.2)(-1.3,-.3)(-1,0)

\pscircle(0,0){1}
\psdot(0,-1)

\pscurve(1,0)(1.2,.2)(1.7,-.3)(2,0)

\pscircle(3,0){1}
\vx{3,-1}{C}

\pscircle(6,0){1}
\psdot(6,1)
\vx{6,-1}{A}

\pscurve(4,0)(4.2,.2)(4.7,-.3)(5,0)

\pscircle(9,0){1}
\rput(9,0){
\psdot(1;300)
\vx{1;60}{B}
\vx{1;180}{D}
}
\Bcedg{A}{D}
\Dcedg{C}{B}



\end{pspicture}
}

\subfigure[\label{f.type3}The planar player has added an edge to produce the Type 3 component.]{
\begin{pspicture}(-1,-1)(10,1)
\SpecialCoor

\pscircle(0,0){1}
\psdot(0,1)
\psdot(0,-1)

\pscurve(1,0)(1.2,.2)(1.7,-.3)(2,0)

\pscircle(3,0){1}
\psdot(3,1)
\vx{3,-1}{A}

\pscircle(6,0){1}
\psdot(6,1)
\vx{6,-1}{B}

\pscurve(7,0)(7.2,.2)(7.7,-.3)(8,0)

\pscircle(9,0){1}
\psdot(9,1)
\psdot(9,-1)

\Dcedg{A}{B}
\end{pspicture}
}

\subfigure[\label{f.type3type3}The planar player has added an edge to the Type 3 component.]{
\begin{pspicture}(-4,-1.55)(7,1.6)
\SpecialCoor

\pscircle(-3,0){1}
\psdot(-3,-1)
\vx{-3,1}{A}

\pscurve(-2,0)(-1.8,.2)(-1.3,-.3)(-1,0)

\pscircle(0,0){1}
\vx{0,-1}{C}

\pscurve(1,0)(1.2,.2)(1.7,-.3)(2,0)

\pscircle(3,0){1}
\psdot(3,1)

\pscircle(6,0){1}
\vx{6,-1}{D}
\vx{6,1}{B}

\pscurve(4,0)(4.2,.2)(4.7,-.3)(5,0)

\Dcedg C D

\end{pspicture}
} \hspace{.5cm}%
\subfigure[\label{f.type3type3a}The planar player has added an edge to the Type 3 component.]{
\begin{pspicture}(-4,-1.55)(4,1.6)
\SpecialCoor

\pscircle(-3,0){1}
\psdot(-3,-1)
\vx{-3,1}{A}

\pscurve(-2,0)(-1.8,.2)(-1.3,-.3)(-1,0)

\pscircle(0,0){1}
\vx{0,-1}{C}

\pscurve(1,0)(1.2,.2)(1.7,-.3)(2,0)

\pscircle(3,0){1}
\vx{3,1}{B}
\psdot(3,-1)

\Dcedg{B}{A}

\end{pspicture}
}

\end{center}
\caption{\label{f.moves} In each case, the nonplanar player has joined two components with the dark edge, and the planar player replies with the dashed edge.}
\end{figure}

Our inductive argument hinges on the fact that when a planar graph has
a bridge, any drawing of it can be ``flipped'' along the bridge to adjust the order of vertices appearing on its outer face.
\begin{observ}
If $G$ is drawn in the plane such that vertices appear in the cyclic order $v_1,\dots,v_k,v_{k+1},\ldots,v_\ell,v_{k+1},v_k,v_1$ along the outer face (so $v_kv_{k+1}$ is a bridge of $G$), then $G$ can also be drawn in the plane such that the vertex order is $v_1,\dots,v_k,v_{k+1},v_\ell,v_{\ell-1},\dots,v_{k+1},v_k,v_1$ along the outer face.  (Note that the $v_i$'s are not necessarily all distinct.)  \qed
\label{o.flip}
\end{observ}

Since all components are of Type \ref{t1} at the beginning of the game, and of Type \ref{t1} or Type \ref{t2} after the first turn of the game, we assume by induction that
condition $\cT$ holds and show that the planar player can respond to any move by the nonplanar player to preserve condition $\cT$.
This will prove Theorem \ref{t.3game}.

We show that the planar player can maintain his invariant via the following remaining cases:
\begin{enumerate}
\item If the nonplanar player has added an edge to a component $C$ of Type 1 or Type 2, the result is already of Type 1; if he has joined two Type 1 components, the result is already of Type 2 or Type 1.  The planar player thus has a free move, which will be addressed in case \ref{p.free}.
\item If the nonplanar player has joined a Type 2 component $C$ to a Type 1 component $C'$, the planar player makes the result a Type 1 component.
\item If the nonplanar player has joined two Type 2 components $C,C'$, the planar player can make the result a Type 2 component.
\item If the nonplanar player has joined the Type 3 component $C$ to a Type 2 component $C'$, the planar player can make the result the Type 3 component.
\item If the nonplanar player has joined the Type 3 component $C$ to a Type 1 component $C'$, the planar player can make the result a Type 3 or Type 2 component.
\item If the nonplanar player has added an edge within the Type 3 component $C$, the planar player can make the result a Type 1 component.
\item \label{p.free} In any other case (a ``free move'' for the planar player) he either turns a Type 3 component into a Type 2 or Type 1 component, or creates a Type 3 component.
\end{enumerate}

\noindent\textbf{Case proofs:}
\nopagebreak
\begin{enumerate}
\item If he adds an edge to a Type 1 or Type 2 component, the deficit of the component after the nonplanar player's move is $\leq 4-2=2$ (e.g., as in \fref{type2self}). If he joins two Type 1 components, the result is a Type 2 or Type 1 component, as in \fref{type1type1}.
\item 
We have that $C$ decomposes as $C_1,C_2$ joined by a bridge, each of whose vertex sets have deficit equal to 2, and that $C$ can be drawn in the plane with all positive deficit vertices on the unbounded face.  Assume that the nonplanar player's move is an edge from $C_2$ to $C'$. Then the planar player chooses an edge between $C_1$ and $C'$, as in \fref{type2type1} (unless $C'$ now has deficit 0, in which case the planar player has a free move).  
The result is a Type 1 component.

\item 
Let the sides of $C$ and $C'$, respectively be $C_1,C_2$ and $C_1'$, $C_2'$.   Suppose without loss of generality that the nonplanar player's move is to take an edge from $C_2$ to $C_1'$.  The planar player responds with a move from $C_2$ to $C_2'$, as shown in \fref{type2type2}.

\item 
Letting $C_i$ ($i=0,1,2$ or $i=0,1,2,3$) be as in the definition of a
Type 3 component, we may assume without loss of generality that the
nonplanar player's edge is from either $C_0$ or $C_1$ to $C'$.  In the
first case, the planar player responds with an edge from $C_1$ to
$C'$, as in \fref{type3type2}.  In the second case, the planar player
takes an edge from $C_0$ to $C'$ (the reverse of the case shown \fref{type3type2}).  Either way, the result is still a single component of Type 3 (of the same subtype \ref{small3} or \ref{big3} as before the nonplanar player's move).  Note that the fact that the resulting component satisfies the condition that all positive-deficit vertices can be drawn on the unbounded face is a consequence of \oref{flip}, which, applied to the Type 3 component $C$, implies that $C$ can be drawn such that the two vertices in $C$ incident with the two new edges are consecutive along the outer face of $C$, among positive deficit vertices.

\item 
This situation is analogous to the previous one.  Again, the first case is shown, in \fref{type3type1}.  The result is either a component of Type \ref{small3} or Type 2; \oref{flip} is used in the same way.

\item If $C$ was of Type \ref{small3}, then it is already a component of Type 1 after the nonplanar player's move, so the planar player has a free move.   If $C$ was of Type \ref{big3}, then it is always possible for the planar player to add a second edge to the component, since, a Type \ref{big3} component admits edges between both the pairs $C_0,C_2$ and $C_1,C_3$ of its components under removal of the its bridges $e_1,e_2,e_3$, and \oref{flip} implies that these edges can be added while preserving the property that the result can be drawn in the plane with all positive-definite vertices on the outer face.   Adding the second edge brings the deficit to 2, making the result again a Type 1 component.

\item
If there remain any two Type 1 components of positive deficit, the planar player can join them to produce a Type 2 (or Type 1) component. 
Otherwise, if there is already a component $C$ of Type 3, he can add
an edge to $C$ to produce a component of Type 1 or Type 2, as in
Figures \ref{f.type3type3} and \ref{f.type3type3a}.  (In the first
case \oref{flip} ensures that the result can be drawn as indicated.)
Or if there is currently no Type 3 component, the planar player can
join a Type 1 or Type 2 component $C$ to a Type 2 component $C'$ to
produce a single component of Type 1, 2, or 3 (depending on the type
of $C$); the case where $C$ has Type 2 is shown in \fref{type3}.
Finally, if no move described so far is possible because there is at
most one component remaining in the graph, and this component is of
Type 1 or Type 2, then he can make arbitrary legal moves until the end
of the game without endangering planarity. This completes the proof of
Theorem \ref{t.3game}.\qed
\end{enumerate}

\bigskip
We turn now to the proof of \tref{4game}.  We call the player with the goal of forcing a $K_\ell$ minor the \emph{minor player}, and the player with the goal of avoiding this the \emph{structure player}.  Our proof has the following two ingredients:
\begin{lemma}
In the course of playing the 4-regular graph game, a player can force the appearance of components of arbitrarily large deficit.
\label{l.largeC}
\end{lemma}

\begin{lemma}
Suppose $G$ is a connected labeled graph, with nonnegative vertex labels bounded some fixed number $b$.  For any $s$, if the sum of the labels of $G$ is sufficiently large relative to $b$ and $\Delta(G)$, we can find $k$ disjoint connected subgraphs of $G$ each with label sums $\geq s$.
\label{l.manyK}
\end{lemma}

By \lref{largeC}, the minor player can build arbitrarily large deficit components in the course of play.  Applying \lref{manyK}, we see that he can find $\ell$ disjoint connected subgraphs each of total deficit $\geq \binom{\ell}{2}$.  Over the next at most $\binom{\ell}{2}$ moves, the minor player joins previously unconnected pairs of these $\ell$ subgraphs (the deficit of each subgraph will remain positive while he is not yet finished), creating a $K_\ell$ minor.\qed

All that remains is to prove Lemmas \ref{l.largeC} and \ref{l.manyK}.
\begin{proof}[Proof of \lref{largeC}]
First note that this Lemma would be very easy if we were instead considering the 5-regular graph game, as then the minor player could simply grow an arbitrarily large deficit component by joining it to isolated vertices on each of his turns; while following such a strategy, the deficit of the component increases by at least $5-2-2=1$ after each time both of the players have made a move.  For the 4-regular graph game, this Lemma will require a bit of care; note, for example, that the Lemma does \emph{not} hold for the 3-regular graph game, even though the invariant the planar player maintains to win that game allows the presence of arbitrarily many deficit-4 components.   (In particular, having isolated vertices of deficit 4 is a stronger condition than having general components of deficit 4.)

We divide the minor player's strategy into two rounds.  In the first round, he chooses $m$ edges of a matching for some large $m$.  Note that if we ignored the role of the structure player, the result would be a large number of components, each of size 2, of deficit 6.  We let $C_1,C_2,\dots,C_m$ denote these pairs of vertices (as sets) which the minor player has joined.  

We claim that once the minor player has completed this round, the sum
\[
\delta= \sum_{\df(C)\geq 5}\df (C)
\]
(taken over all connected components of the graph which at this round have deficit $\geq 5$) is large (tends to $\infty$ with $m$).   To see this, let us allow even that the structure player is given the $m$ edges of the minor player's matching ahead of time.  In absence of the structure player's moves, $\delta$ would be $6m$.  We classify the structure player's moves in this round into two types:
\begin{enumerate}
\item Edges in components which, at the end of this round have deficit $\leq 4$, and
\item Edges in components which, at the end of this round have deficit $\geq 5$.
\end{enumerate}
Let $m_i$ denote the number of moves he makes of Type $i$ (so $m=m_1+m_2$).   We have that
\[
\delta\geq 6m-6m_1/\beta-2m_2,
\]
where the constant $\beta$ is the minimum number of edges per component required to decrease the deficit of components below 5.  We need only show that $\beta>1$.  To bound $\beta$, consider any component $C$ of the graph after this first round of play.  If it contains $m_C$ edges of the minor player's matching and $m_C'$ edges of the structure player, then its deficit satisfies
\[
\df(C)\geq 6m_C-2m_C'.
\]
In particular, $\df(C)\leq 4$ implies that
\[
m'_C\geq 3m_C-2.
\]
Moreover, in the case where $m_C=1$, we see that $m'_C\geq 5$ since the structure player cannot duplicate the minor player's edge (his best case is to complete a $K_4$).  In particular, since $m_C'\geq 4$ in all cases, $\beta\geq 4$, completing the proof that $\delta$ becomes arbitrarily large.

\smallskip
We now show that he can force the appearance of a single component of large deficit.  He simply chooses one component of deficit $\geq 5$ arbitrarily, and, on each turn, grows this component by joining it to a new component of deficit $\geq 5$ arbitrarily.  Taking into account also the structure player's move, the deficit of this component is increasing by at least $5-2-2>1$ on each turn; in particular, it will become arbitrarily large.
\end{proof}
\begin{proof}[Proof of \lref{manyK}]
Let $G$ be a connected graph with maximum degree $\Delta$ with vertices with labels $\ell(v)$ from $0,1,\dots,b$, and let $\ell(X)$ denote the label sum of a subset $X\sbs G$.  Consider a spanning tree $T$ of $G$.  We begin by showing that when $\ell(G)$ is sufficiently large, we can find an edge $e$ of $T$ such that the ratio $\rho_e=\ell(C_2)/\ell(C_1)$ for the components $C_1,C_2$ of $T\stm e$ satisfies
\begin{equation}
\frac 1 \Delta < \rho_e < \Delta.
\label{e.rhorange}
\end{equation}
To see this, fix an edge $e\in T$, and let us consider the case where $\rho_e<\frac 1 \Delta$.  If $x$ is the endpoint of $e$ in the larger label-sum component $C_1$ of $T\stm e$, and $e'=\{x,y\}$ is the edge for which $y$ is in the highest label-sum component of $C_1\stm x$, then we have that
\[
\rho_{e'}\leq \frac{\ell(C_2)+b+\frac{\Delta-2}{\Delta-1}\ell(C_1)}{\frac 1 {\Delta-1}(\ell(C_1)-b)}
\]
Letting $\ell(G)$ be sufficiently large that, say, $b<(\frac 1 {\Delta-1}-\frac 1 \Delta)\ell(C_1)$ (so that also $b<\frac {\ell(C_1)} \Delta$) 
we see that
\[
\rho_{e'}\leq \frac{\Delta}{\Delta-1}\frac{\ell(C_2)+\frac{\Delta-1}{\Delta}\ell(C_1)}{\frac{1}{(\Delta-1)}\ell(C_1)}=(\Delta-1)+\Delta\rho_e.
\]
In particular, if $\rho_e<\frac 1 \Delta$, then $\rho_e<\rho_{e'}<\Delta$; thus, we can walk along the tree to find an edge $e$ satisfying \eref{rhorange}.

We now simply apply our ability to find such edges recursively, $t$ times for some $t$, to divide $T$ into $2^t$ trees, each with label sums $\geq \frac {\ell(G)} {(\Delta+1)^t}.$
\end{proof}
\section{Discussion}
Although we have focused on topological questions regarding the game we have introduced, many other questions seem natural as well.  For example, what if we consider subgraphs instead of minors?  For example:
\begin{q}
Which graphs $H$ have the property that a player in the 3-regular graph game on sufficiently many vertices can ensure that the resulting graph contains a copy of $H$?
\end{q}
For example, $K_4$ does not have this property; indeed no cubic graph has this property.  This can be seen by showing that either player in the 3-regular graph game has a strategy to ensure that the graph is connected.  (The player can maintain the invariant that while disconnected, the graph consists of isolated vertices, plus a single other component with at least one vertex of deficit 2.)

Turning back to topological issues, one can probe the relationship between degree-4 and degree-3 vertices a bit more.  Let us define the \emph{cubic$+k$ graph game}, in which degrees of vertices must remain at most 3, except for $k$ special vertices whose degrees may rise to 4.  Let now $g^3(k)$ be the minimum genus $g$ such that for any number of vertices $n$, Player 1 can ensure that any connected component of the result of the cubic$+k$ game can be drawn in some surface of genus $g$.   Then Theorem \ref{t.3game} implies that $g^3(0)=0$, while the proof of Theorem \ref{t.4game}, which works when some vertices have a degree threshold of 3 so long as sufficiently many have a degree threshold of 4, implies that $g^3(k)\to \infty$.  Thus the asymptotic behavior of $g^3(k)$ can be studied.  Rather than the particular rate of growth of $g^3(k)$, however, it may be more interesting to compare with the function $g_4(n)$, which we define as the minimum genus such that Player 1 can ensure that the any connected component of the result of the 4-regular graph game on $n$ vertices can be drawn on some surface of genus $g$:
\begin{q}
Is $g^3(n)\sim g_4(n)$?
\end{q}
\noindent A first step would be finding a single value of $n$ for which $g^3(n)\neq g_4(n)$.


\begin{thebibliography}{99}
\bibitem{archdeacon1981kuratowski} D. Archdeacon.  A Kuratowski theorem for the projective plane, in {\em Journal of Graph Theory} {\bf 5} (1981) 243--246.


\bibitem{Beck} J. Beck, Combinatorial Games: Tic-Tac-Toe Theory,
  Cambridge University Press, 2008.

\bibitem{glover1975cubic} H. Glover and J. Huneke.  Cubic irreducible graphs for the projective plane, in {\em Discrete Mathematics} {\bf 13} (1975) 341--344

\bibitem{grundy} P. Grundy, Mathematics and games, in {\em Eureka} {\bf 2} (1939) 21.
\bibitem{sprague} R. Sprague, Uber mathematische kampfspiele, in {\em  T{\^o}hoku Math. J} {\bf 41} (1935) 438--444.






\end{thebibliography}
\end{document}